\documentclass[12pt]{amsart}
\usepackage{amscd}
\usepackage{graphicx}
\def\frk{\frak}               

\def\Phi{{\frk n}}
\def\Phi{{\frk N}}
\def\opn#1#2{\def#1{\operatorname{#2}}} 
\opn\chara{char} \opn\length{\ell} \opn\pd{pd} \opn\rk{rk}
\opn\projdim{proj\,dim} \opn\injdim{inj\,dim} \opn\rank{rank}
\opn\depth{depth} \opn\grade{grade} \opn\height{height}
\opn\embdim{emb\,dim} \opn\codim{codim}

\opn\Tr{Tr} \opn\bigrank{big\,rank}
\opn\superheight{superheight}\opn\lcm{lcm}
\opn\trdeg{tr\,deg}
\opn\reg{reg} \opn\lreg{lreg} \opn\ini{in} \opn\lpd{lpd}
\opn\size{size}
\opn\div{div} \opn\Div{Div} \opn\cl{cl} \opn\Cl{Cl}
\opn\Spec{Spec} \opn\Supp{Supp} \opn\supp{supp} \opn\Sing{Sing}
\opn\Ass{Ass} \opn\Min{Min}
\opn\Ann{Ann} \opn\Rad{Rad} \opn\Soc{Soc}
\opn\Im{Im} \opn\Ker{Ker} \opn\Coker{Coker} \opn\Am{Am}
\opn\Hom{Hom} \opn\Tor{Tor} \opn\Ext{Ext} \opn\End{End}
\opn\Aut{Aut} \opn\id{id}

\opn\nat{nat}
\opn\pff{pf}
\opn\Pf{Pf} \opn\GL{GL} \opn\SL{SL} \opn\mod{mod} \opn\ord{ord}
\opn\Gin{Gin} \opn\Hilb{Hilb}
\opn\aff{aff} \opn\con{conv} \opn\relint{relint} \opn\st{st}
\opn\lk{lk} \opn\cn{cn} \opn\core{core} \opn\vol{vol}
\opn\link{link} \opn\star{star}
\opn\gr{gr}

\def\pot#1#2{#1[\kern-0.28ex[#2]\kern-0.28ex]}

\opn\dirlim{\underrightarrow{\lim}}
\opn\inivlim{\underleftarrow{\lim}}
\def\Implies{\ifmmode\Longrightarrow \else
        \unskip${}\Longrightarrow{}$\ignorespaces\fi}
\def\implies{\ifmmode\Rightarrow \else
        \unskip${}\Rightarrow{}$\ignorespaces\fi}
\def\iff{\ifmmode\Longleftrightarrow \else
        \unskip${}\Longleftrightarrow{}$\ignorespaces\fi}

\let\:=\colon
\newtheorem{Theorem}{Theorem}[section]
\newtheorem{Lemma}[Theorem]{Lemma}
\newtheorem{Corollary}[Theorem]{Corollary}
\newtheorem{Proposition}[Theorem]{Proposition}
\newtheorem{Remark}[Theorem]{Remark}

\newtheorem{Example}[Theorem]{Example}

\newtheorem{Definition}[Theorem]{Definition}

\let\epsilon\varepsilon
\let\phi=\varphi
\let\kappa=\varkappa
\def\qed{\ifhmode\textqed\fi
      \ifmmode\ifinner\quad\qedsymbol\else\dispqed\fi\fi}
\def\textqed{\unskip\nobreak\penalty50
       \hskip2em\hbox{}\nobreak\hfil\qedsymbol
       \parfillskip=0pt \finalhyphendemerits=0}
\def\dispqed{\rlap{\qquad\qedsymbol}}

\opn\dis{dis}
\def\pnt{{\raise0.5mm\hbox{\large\bf.}}}

\opn\Lex{Lex}

\begin{document}
\title{Inclusion Ideals associated to Uniformly increasing Hypergraphs}
 \maketitle

\author\begin{center}{I. Anwar$^{a}$, S. Ahmad$^{a}$, A. Inam, A. Haider$^{b}$}\end{center}
\begin{center}\footnotesize{\textsl{{
$^{a}$ COMSATS Institute of Information Technology Lahore, Pakistan.\\ E. mails: imrananwar@ciitlahore.edu.pk, sarfrazahmad@ciitlahore.edu.pk, amina.inam@gmail.com\\
$^{b}$ JAZAN University, Saudi Arabia.\\ E.
mail:azeemhaider@gmail.com}}}
\end{center}

\begin{abstract}
In this paper, we introduce {\em inclusion ideals} $\mathcal
I(\mathcal H)$ associated to a special class of {\em non uniform
hypergraphs $\mathcal{H}(\mathcal{X},\mathcal{E},d)$, }namely, the
{\em uniformly increasing hypergraphs}. We discuss some algebraic
properties of the inclusion ideals. In particular, we give an upper
bound of the {\em Castelnouvo-Mumford regularity} of the {\em
special dual ideal $\mathcal I^{[*]}(\mathcal H).$}
 \vskip 0.4 true cm
 \noindent
\end{abstract}

\section{introduction}

Let $\mathcal{E}=\{{E_{1},...,E_{s}}\}$ be the collection of
distinct subsets of the finite set
$\mathcal{X}=\{{x_{1},...,x_{n}}\}$. The pair
$\mathcal{H}(\mathcal{X},\mathcal{E})$ is said to be a {\em
Hypergraph}, if $E_{i}\neq \emptyset $ for each $i,$ where
$\mathcal{X}$ and $\mathcal{E}$ are the set of {\em vertices} and
{\em edges} of $\mathcal{H}$, respectively. A hypergraph is said to
be {\em d-uniform hypergraph}, if $|E_{i}|=d$ for each $i$. For
example, all simple graphs are {\em 2-uniform hypergraphs}. It is
worth mentioning here that hypergraphs are the
generalized form of the simple graphs and the simplicial complexes.\\
Let $S=K[x_{1},...,x_{n}]$ be the polynomial ring over an infinite
field $K$. If we associate each vertex $x_{i}$ to each variable
$x_{i}$, then the edge ideal $ I(\mathcal H)$ associated to simple
hypergraph (see \cite{Ha}) $\mathcal{H}$ is:
$$ I(\mathcal H)=(\{x^E=\prod_{x\in E} x \mid E\in \mathcal{E}\})\subseteq S$$
This edge ideal was introduced by  Villarreal for a special case of
hypergraphs, namely, the simple graphs. The edge ideal $ I(\mathcal
H)$ of a hypergraph was firstly discussed by Faridi in \cite{F1} and
\cite{F2}, but in different context. The study of algebraic objects
via combinatorial correspondence has fascinated many people, who
 have been working to build a dictionary between the
algebraic properties of edge ideal and the combinatorial structure
associated to it. For instance, the description about the edge
ideals and the ideals associated to simplicial complexes can be
found in \cite{F1}, \cite{F2}, \cite{Mi}, \cite{St},
 \cite{St1} and \cite{Vi}. \\
In this paper, we introduce the monomial ideal $\mathcal I(\mathcal
H)$ associated a special kind of non uniform Hypergraphs namely {\em
uniformly increasing hypergraphs}
$\mathcal{H}(\mathcal{X},\mathcal{E},d)$, see \ref{uni}. We call
$\mathcal I(\mathcal H)$ as the {\em inclusion ideal}, which is not
a square-free ideal in contrast to the edge ideal, see \ref{inc}.
 In \ref{TO}, we prove that the {\em special dual} of inclusion ideal
$\mathcal{I}^{[*]}(\mathcal{H})$ is the monomial ideal whose
associated prime ideals are totally ordered by inclusion. As can be
easily seen, by an appropriate change of variables, such ideals are
the monomial ideals of Borel type introduced by Herzog, Popescu and
Vladoiu in \cite{HPV}. Moreover, Herzog and Popescu in \cite{HP}
proved that such ideals are {\em pretty clean}. The upper bound for
the {\it Castlenouvo-Mumford regularity} of such ideals are
discussed in \cite{AA} and \cite{Ci}. The  {\it Castlenouvo-Mumford
regularity} of $I$ is given by $reg(I)=max\{j-i|\beta_{ij}(I)\neq
0\}$, if $\beta_{ij}(I)$ are the graded Betti numbers of $I$.
  In \ref{Reg} we
give the upper bound for the regularity of
$\mathcal{I}^{[*]}(\mathcal{H})$, which is more finer than the bound
already found by Ahmad, Anwar in \cite{AA} and Cimpoeas in
\cite{Ci}.\\
We would like to thank the referee for his valuable suggestions and
remarks.
\section{Inclusion Ideal and its special dual}
This section is devoted for the introduction of {\em uniformly
increasing hypergraphs} $\mathcal{H}(\mathcal{X},\mathcal{E},d)$ and
its inclusion ideal we associate to it.
\begin{Definition}\label{uni}
{\em Let $\mathcal{X}=\{{x_{1},...,x_{n}}\}$ be a finite set and
$\mathcal{E}=\{{E_{1},...,E_{s}}\}$ be a collection of subsets of
$\mathcal{X}$ such that $E_{i}\subseteq E_{j}$ for $1\leq i< j\leq
s$, where $E_{i},E_{j}\in \mathcal{E}$, then the triplet
$\mathcal{H}(\mathcal{X},\mathcal{E},d)$ is called
 {\em uniformly increasing hypergraph}, if\\
(1) $|E_{i}|\geq 2$ and\\
(2)$|E_{i+1}|=|E_{i}|+d$, where $d\in \mathbb Z^+$ is the
increment.}
\end{Definition}
\begin{Example}\label{ex1}
Let $\mathcal X=\{x_1,x_2,x_3,x_4\}$ be the set of vertices and
$\mathcal E=\{E_1,E_2,E_3\}$ be the set of edges, where
$E_1=\{x_1,x_2\}, \, E_2=\{x_1,x_2,x_3\}, \,
E_3=\{x_1,x_2,x_3,x_4\}$, where $|E_{i+1}|-|E_{i}|=1$ for $1\leq i
\leq 3$, then the hypergraph $\mathcal H(\mathcal X,\mathcal E,1)$
is:
\begin{figure}[h]
\begin{center}
  \includegraphics[width=6.2cm]{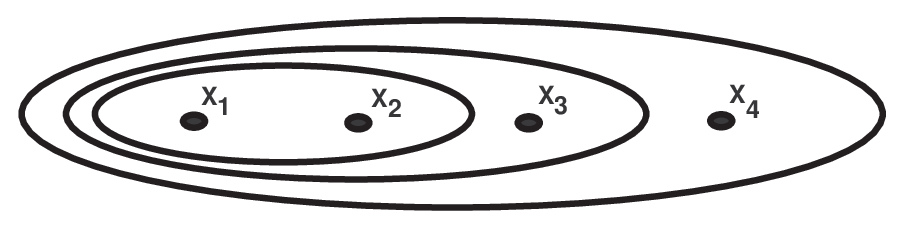}\vspace{0.0cm}
  \label{1}
\end{center}
\end{figure}
\end{Example}
\begin{Definition}\label{co}{\em
For a {\em uniformly increasing hypergraph}
$\mathcal{H}(\mathcal{X},\mathcal{E},d)$ over $n$-vertices, we
define the {\em containment vector} `$a$' as $a=(a_1, a_2, \ldots,
a_n)$ where each $a_i$ is the {\em degree} of vertex $x_i$; that is
$$a_i=\mid \{E_j \, \mid \, x_i\in E_j \hbox{ with }
j\in\{1,2,\ldots,s\} \} \mid $$ where $\mid \, . \mid$ denotes the
cardinality of the set.

}\end{Definition}
\begin{Example}\label{A1}{\em

 For the hypergraph in example \ref{ex1}; the {\em containment vector} is\linebreak $a=(3,3,2,1)$. As  $x_1$
belongs to $\{E_1, E_2, E_3\}$, so $a_1=3$. Similarly, we have
$a_2=3$, $a_3=2$ and $a_4=1.$
 }
\end{Example}
\begin{Definition}\label{inc}
\em{Let $S=k[x_{1},...,x_{n}]$ be a polynomial ring over an infinite
field $k$ and let $\mathcal H(\mathcal X,\mathcal E,d)$ be a {\em
uniformly increasing hypergraph} over $n$ vertices. Considering the
correspondence of each vertex $x_{i}$ of $\mathcal H(\mathcal
X,\mathcal E,d)$ to each variable $x_{i}$ in $S$, the {\em inclusion
ideal} $\mathcal I(\mathcal H)$ of the {\em uniformly increasing
hypergraph} is defined as,
 $$\mathcal I(\mathcal H)=(\{x^{E_i}=\prod_{x\in E_i} x^{s-i+1} |E_i\in \mathcal E,\,\,
\forall,\,\,1\leq i\leq s=|\mathcal{E}| \}).$$ It should be noted
that $\mathcal I(\mathcal H)$ is minimally generated by the
above\linebreak mentioned monomial generators and we will denote it
with $G(\mathcal I(\mathcal H))$. }\end{Definition}

\begin{Example}\label{a}
 {\em Let $\mathcal H(\mathcal X,\mathcal E,d)$ be the hypergraph in
 above example \ref{ex1}, Then by definition its {\em inclusion ideal} $\mathcal I(\mathcal H)$
 will be:
$$\mathcal I(\mathcal H)=(x_1^3x_2^3,x_1^2x_2^2x_3^2,x_1x_2x_3x_4).$$
}
\end{Example}
 \begin{Remark}{\em
Note that, the {\em inclusion ideal} $\mathcal I(\mathcal H)$ of a
{\em uniformly increasing hypergraph} $\mathcal H(\mathcal
X,\mathcal E,d)$ on $n$ vertices is a non square-free monomial ideal
in $S$ unlike the edge ideals associated to the uniform hypergraphs.
That is why, we call the {\em inclusion ideal} instead of edge
ideal.}
\end{Remark}

Now we give the definition of the {\em Alexander dual} of an
arbitrary monomial ideal $I$ with respect to `$a$' given in
(\cite{Mi}, section $5.2$).
\begin{Definition}\label{Al} {\em
Given two vectors $c, d\in \mathbb{N}^{n}$ with $d\preceq c$ (that
is $d_i\leq c_i$ for $i=1,\ldots, n$), let $c \setminus d$ denote
the vector in $\mathbb{N}^{n}$whose $i^{th}$ coordinate is
$$c_i\setminus d_i=\left\{
                     \begin{array}{ll}
                       c_i +1-d_i, & \hbox{if } d_i\geq 1; \\
                       0, & \hbox{if  } d_i= 0.
                     \end{array}
                   \right.
$$
If $I\subset S$  is a monomial ideal whose minimal generators all
divide $x^c$, then the {\it {\bf Alexander dual} }of $I$ with
respect to {\it {\bf c}} is
$$I^{[c]}=\bigcap (m^{c\setminus b} |\, \, x^b\,\, \hbox{is the minimal generator of
I}),$$ where each $m^{c\setminus b}$ is the irreducible ideal
generated by powers of variables. That is for a vector ${c\setminus
b}$ in $\mathbb{N}^n$, we have
$$m^{c\setminus b}=<x_i^{c_i\setminus b_i}\mid c_i\setminus b_i\geq 1>$$}
\end{Definition}
\begin{Lemma}{\em
If  $\mathcal I(\mathcal H)=(x^{b_1},x^{b_2}, \ldots , x^{b_s})$
with $b_i\in \mathbb{N}^{n}$ is the inclusion ideal of a {\em
uniformly increasing hypergraph} $\mathcal H(\mathcal X,\mathcal
E,d)$ on $n$ vertices , then $b_i\preceq a$ for all
$i\in\{1,\ldots,s \}$, where `$a$' is the {\em containment vector }
associated to $\mathcal H(\mathcal X,\mathcal E,d)$ . }\end{Lemma}
\begin{proof}
 Let $a=(a_1,\ldots,a_n)$ be the {\em containment vector} associated to
$\mathcal H(\mathcal X,\mathcal E,d)$ such that each $a_i=deg(x_i)$
follows from \ref{co}. If $x_k\in E_t$ such that $x_k\not\in
E_{t-1}$ for some $1\leq t\leq s$, then $x_k\in E_{t+1},\, x_k\in
E_{t+2},\ldots, x_k\in E_{s}$ which implies  $a_k=deg(x_k)=s-t+1$.
Let us take an arbitrary vector $b_i=(b_{i1},b_{i2},\ldots,b_{in})$
in $\mathbb{N}^n$ such that $x^{b_i}\in G(\mathcal I(\mathcal H))$.
Then from \ref{inc}, we have $b_{ij}=0$ if $x_j\not\in E_i$ or
$b_{ij}=s-i+1$ if $x_j\in E_i$ for all $j\in\{1,\ldots,n\}$. If
$b_{ij}=0$, then $b_{ij}<a_j$ because $a_j>0$. If $b_{ij}=s-i+1$,
then\linebreak $ s-i+1\leq s-t+1= a_j$, for $t\leq i$ such that
$x_j\in E_t$ with $x_j\not\in E_{t-1}$.
\end{proof}

Next we give the definition of the {\em special dual} of the
inclusion ideal of a {\em uniformly increasing hypergraph}.
\begin{Definition}\label{sp}{\em
If  $\mathcal I(\mathcal H)=(x^{b_1},x^{b_2}, \ldots , x^{b_s})$ is
the inclusion ideal of a {\em uniformly increasing hypergraph}
$\mathcal H(\mathcal X,\mathcal E,d)$ on $n$ vertices, then the
Alexander dual of $\mathcal I(\mathcal H)$ with respect to the {\em
containment vector} `$a$' associated to the hypergraph $\mathcal
H(\mathcal X,\mathcal E,d)$ is called the {\em special dual} of the
inclusion ideal, and is denoted by $I^{[*]}(\mathcal H)$. That is;
$$I^{[*]}(\mathcal H)=\bigcap_{i=1}^{s} (m^{a\setminus b_i} |\, \, x^{b_i}\,\, \hbox{is the minimal generator of
I}).$$ It should be noted that here $x^{b_i}$ is actually $x^{E_i}$.
  }\end{Definition}

\begin{Example}
Let $\mathcal I(\mathcal H)$ be the monomial ideal in above example
\ref{a}, its {\em special dual} will be:
 $$\mathcal I^{[*]}(\mathcal H)=(x_1,x_2)\cap (x_1^2,x_2^2,x_3)\cap
(x_1^3,x_2^3,x_3^2,x_4)$$
$$=(x_1^3,x_2^3,x_1x_3^2,x_2x_3^2,x_1^2x_4,x_2^2x_4,x_1x_3x_4,x_2x_3x_4)$$
here {\em containment vector} is $a=(3,3,2,1)$.

\end{Example}

For a monomial $m\in S=k[x_1, x_2,\ldots,x_n]$, we define $\supp(m)$
as\linebreak $\supp(m)=\{x_i \mid \, x_i\mid m \hbox{ \, for\,  }
i\in\{1,2,\ldots,n\} \}.$ Also for a monomial ideal $I=(g_1,\ldots
,g_m)\subset S$ with $\{g_1,\ldots, g_m\}$ is the minimal set of
generators of $I$, we define
$\supp(I)=\cup_{i=1}^{m}\{\supp(g_i)\}$.

\begin{Proposition}\label{TO}
{\em Let $\mathcal I(\mathcal H)$ be an inclusion ideal in
$S=k[x_1,...,x_n]$ associated to a {\em uniformly increasing
hypergraph} $\mathcal H(\mathcal X,\mathcal E,d)$ on $n$ vertices,
then $Ass(S/\mathcal I^{[*]}(\mathcal H))$ is totally ordered by
inclusion. }
\end{Proposition}

\begin{proof}

 Let $\mathcal I(\mathcal H)=(x^{b_1}, x^{b_2},\ldots, x^{b_s}) \in S$ be the inclusion ideal associated to a uniformly
 increasing hypergraph $\mathcal H(\mathcal X,\mathcal E,d)$ on $n$ vertices. So, the irredundant irreducible primary
 decomposition of special dual $\mathcal I^{[*]}(\mathcal H)$ (from the definition
\ref{sp}) is as follows;
  $$\mathcal I^{[*]}(\mathcal H)=\bigcap (m^{a\setminus b_i} |\, x^{b_i}\,\, \hbox{ be\,\, the\,\,
  minimal\,\,
generators\,\, of\,\, } \mathcal I).$$ For each $x^{b_i}\in
G(\mathcal I(\mathcal H))$, we have an irreducible primary ideal
$m^{a\setminus b_i}$ appearing in the above primary decomposition of
$\mathcal I^{[*]}(\mathcal H)$ with
$\supp(x^{b_i})=\supp(m^{a\setminus b_i})$. As each generator
$x^{b_i}$ of the inclusion ideal $\mathcal I(\mathcal H)$ is
associated to each edge $E_i$ of the increasing hypergraph. So, we
have $\supp(x^{b_1})\subset \supp(x^{b_2})\subset \ldots \subset
\supp(x^{b_s})$ \linebreak i.e,
 $\sqrt{{m^{a\setminus b_i}}}=\supp(m^{a\setminus b_i})= \supp(x^{b_i})$. Hence the
associated primes of the above primary decomposition of $\mathcal
I^{[*]}(\mathcal H)$ are totally ordered under inclusion.
\end{proof}
\section{Regularity of the special dual $\mathcal I^{[*]}(\mathcal H)$ }
 Let $k$ be an infinite field,
$S=k[x_1,...,x_n],n\geq 2$ the polynomial ring over $k$ and
$I\subset S$ a monomial ideal. Let $G(I)$ be the minimal set of
monomial generators of $I$ and $deg(I)$ the highest degree of a
monomial of $G(I)$. Given a monomial $u\in S$ set $m(u)=max\{i\mid
x_i|u\}$ and $m(I)=max_{u\in G(I)}\ m(u)$.  Also, $I_{\geq t}$ be
the ideal generated by the monomials of $I$ of degree $\geq t$. A
monomial ideal $I$ is {\it stable} if for each monomial $u\in I$ and
$1\leq j<m(u)$ it follows $\frac{x_ju}{x_{m(u)}}\in I$.
\begin{Proposition}\label{L1}{\em
A monomial ideal $I=(x_{1}^{a_1},...,x_{n}^{a_n})\subset S$, with
\linebreak $n-1\geq a_1\geq a_2\geq \ldots a_n\geq 1$ has $I_{\geq
t(I)}$ stable, where $t(I)=\sum_{i=2}^n a_i$. }\end{Proposition}
\begin{proof}
If $u\in I_{\geq t(I)}$, from abov we get $u=v\cdot x_j^{a_j}$ for
some $1\leq j\leq n$ and \linebreak $v\in
(x_1,\ldots,x_n)^{t(I)-a_j}$, then $u$ belongs to the stable ideal
$(x_1,\ldots,x_n)^{t(I)}$ and it is enough to show that
$(x_1,\ldots,x_n)^{t(I)}\subseteq I_{\geq t(I)}$. If $w\in
(x_1,\ldots,x_n)^{t(I)}$, then $w=x_1^{\alpha_1}x_2^{\alpha_2}\ldots
x_n^{\alpha_n}$ with all $\alpha_i\geq 0$ and $\Sigma_{i=1}^{n}
\alpha_i \geq t(I)$. Now we will prove that there exists some $k$
such that $1\leq k\leq n$ with $\alpha_k\geq a_k$. Suppose contrary
that there does not exist such $k$, that is $\alpha_i<a_i$ for all
$i\in \{1,\ldots,n\}$. Because $\Sigma_{i=1}^{n} \alpha_i \geq
t(I)=\Sigma_{i=2}^{n} a_i$, therefore
$$\alpha_1\geq(a_2-\alpha_2)+(a_3-\alpha_3)+\ldots +(a_n-\alpha_n)$$
So we have $\alpha_1\geq n-1$, hence $a_1> n-1$. which is a
contradiction. So we can take above
$v=x_1^{\alpha_1}...x_k^{\alpha_k-a_k}...x_n^{\alpha_n}$ which shows
that $w\in I_{\geq t(I)} $.

\end{proof}
\begin{Remark}\label{R}{\em
In general one cannot get $I_{\geq t(I)-1}$ stable, when
$I=(x_1^{a_1},x_2^{a_2},...,x_n^{a_n})$ with $n-1\geq a_1\geq
a_2\geq \ldots a_n\geq1$. For example, if $n=3$ and
$I=(x_1^2,x_2^2,x_3)$ then $t(I)=3$ and clearly $I_{\geq 2}$ is not
stable.}
\end{Remark}

Now we recall the following result from \cite{AA};

\begin{Proposition}\label{L2}\cite{AA}
{\em If $I,J$ are monomial ideals such that $I_{\geq q(I)}$ and
$J_{\geq q(J)}$ are stable ideals, then $(I\cap J)_{\geq max (q(I)
,\,  q(J))}$ is stable.}
\end{Proposition}
\begin{Lemma}\label{L3}{\em
If $\mathcal I(\mathcal H) \in S$ be the inclusion ideal associated
to a uniformly
 increasing hypergraph $\mathcal H(\mathcal X,\mathcal E,d)$ on $n$ vertices, then $\mathcal I^{[*]}_{\geq{t(I)}}(\mathcal
 H)$ is stable, where $t(I)=\sum_{i=2}^n a_i$ with $a=(a_1,a_2,\ldots,a_n)$ be the containment vector associated to $\mathcal H(\mathcal X,\mathcal E,d)$.
}\end{Lemma}

\begin{proof}
Let $\mathcal I(\mathcal H)=(x^{b_1}, x^{b_2},\ldots,x^{b_s}) \in S$
with $b_i\in \mathbb{N}^n$ be the inclusion ideal associated to a
uniformly increasing hypergraph $\mathcal H(\mathcal X,\mathcal
E,d)$ on $n$ vertices. Then its special dual $\mathcal
I^{[*]}(\mathcal H)$ will be;
  $$\mathcal I^{[*]}(\mathcal H)=\bigcap_{i=1}^{s=\mid\mathcal E \mid} (m^{a\setminus b_i} |x^{b_i}\,\, \hbox{ be\,\, the\,\,
  minimal\,\, generators\,\, of\,\, } \mathcal I(\mathcal H))$$ or
$$\mathcal I^{[*]}(\mathcal H)=\bigcap_{i=1}^{s} Q_i \, \, \, \hbox{ where }\, Q_i=(m^{a\setminus b_i})\, \hbox{\, for all\, }i\in \{1,\ldots, s\}$$
Because $\Ass(S/\mathcal I^{[*]}(\mathcal H))$ totally ordered from
\ref{TO}, so we have
$$\sqrt{Q_1}\subset \sqrt{Q_2}\subset \ldots \subset \sqrt{Q_s} .$$
Also $m^{a\setminus b_i}$ is an irreducible ideal with
$\supp(m^{a\setminus b_i})=\supp(x^{b_i})$ for each \linebreak
$i\in\{1,\ldots,s\}$. So, we have $Q_i
=(x_1^{\alpha_1},x_2^{\alpha_2},\ldots,x_{r_i}^{\alpha_{r_i}})$ with
$a\setminus b_i=(\alpha_1, \alpha_2, \ldots, \alpha_n)$ where
$\alpha_k=0$ for $r_i+1\leq k\leq n$  and $r_i-1\geq \alpha_1=i\geq
\ldots \geq\alpha_{r_i}= 1$ follows from the definitions \ref{uni}
and \ref{sp}. Therefore,
$Q_s=(x_1^{\alpha_1},x_2^{\alpha_2},\ldots,x_n^{\alpha_n})$ with
$n-1\geq \alpha_1=s\geq \alpha_2\geq \ldots\geq \alpha_n=1$. Hence
from the above Propositions  \ref{L1} and \ref{L2}, it immediately
follows that $\mathcal I^{[*]}_{\geq t(I)}(\mathcal H)$ is stable.
\end{proof}

Next we recall a Proposition from \cite{ERT}.
\begin{Proposition}\label{ert}{\em Let $I$ be a monomial ideal and $e\geq deg(I)$ an
integer such that $I_{\geq e}$ is stable. Then $reg(I)\leq e$.}
\end{Proposition}

\begin{Corollary}\label{Reg}{\em
Let $\mathcal I(\mathcal H) \in S$ be the inclusion ideal associated
to a uniformly  increasing hypergraph $\mathcal H(\mathcal
X,\mathcal E,d)$ on $n$ vertices. Then for its special dual ideal
$reg(\mathcal I^{[*]}(\mathcal
 H))\leq t(\mathcal I^{[*]}(\mathcal
 H))=\sum_{i=2}^s a_i$.}
\end{Corollary}
\begin{proof}By the previous lemma \ref{L3} we have $\mathcal I^{[*]}_{\geq t(\mathcal I^{[*]}(\mathcal
 H))}(\mathcal
 H)$ stable, as\linebreak $deg(\mathcal I^{[*]}(\mathcal
 H))= max\{a_1,a_2,\ldots,a_n\}$ and clearly  $t(\mathcal I^{[*]}(\mathcal
 H))\geq deg(\mathcal I^{[*]}(\mathcal
 H))$. Hence we get $reg(\mathcal I^{[*]}(\mathcal
 H))\leq t(\mathcal I^{[*]}(\mathcal
 H))=\sum_{i=2}^s a_i$ by proposition \ref{ert}.
\end{proof}
\begin{Remark}{\em

As $\mathcal I^{[*]}(\mathcal H)$ is a monomial ideal whose
associated prime ideals are totally ordered under inclusion. In
\cite{AA} and \cite{Ci}, the authors have given the bound for the
regularity of such ideals that is $reg(\mathcal I^{[*]}(\mathcal
 H))< q(\mathcal I^{[*]}(\mathcal
 H))=m(\mathcal I^{[*]}(\mathcal
 H))(deg(\mathcal I^{[*]}(\mathcal
 H))-1)+1$.
Therefore, we have $reg(\mathcal I^{[*]}(\mathcal
 H))\leq q(\mathcal I^{[*]}(\mathcal
 H))$. But it is worth noting that $t(\mathcal I^{[*]}(\mathcal
 H))\leq q(\mathcal I^{[*]}(\mathcal
 H))$,
 for instance, for the ideal $\mathcal I^{[*]}(\mathcal H)=(x_1,x_2)\cap(x_1^2,x_2^2,x_3),$ $t(\mathcal I^{[*]}(\mathcal
 H))=3<q(\mathcal I^{[*]}(\mathcal H))=4$.
  Moreover, it should be noted that
the bound for the regularity of $\mathcal{I}^{[*]}(\mathcal{H})$
found in \ref{Reg} can not be applied to any ideal whose associated
prime ideals are totally ordered by inclusion. It is also important
to mention here that $t(\mathcal I^{[*]}(\mathcal
 H))=\sum_{i=2}^{n}a_i$ is a combinatorial term
depends on the {\em containment vector `$a$'} associated to the {\em
uniformly increasing hypergraph} $\mathcal H(\mathcal X,\mathcal
E,d)$ on $n$ vertices. }
\end{Remark}

\end{document}